\documentclass[12pt]{amsart}

\usepackage{psfrag}
\usepackage{color}
\usepackage{tikz}
\usetikzlibrary{matrix,arrows}
\usepackage{graphicx,graphics}
\usepackage{fullpage,amssymb,amsfonts,amsmath,amstext,amsthm,amscd,verbatim,enumerate}
\usepackage[T1]{fontenc}

\begin{document}

\newtheorem{theorem}{Theorem}[section]
\newtheorem{result}[theorem]{Result}
\newtheorem{fact}[theorem]{Fact}
\newtheorem{example}[theorem]{Example}
\newtheorem{conjecture}[theorem]{Conjecture}
\newtheorem{lemma}[theorem]{Lemma}
\newtheorem{proposition}[theorem]{Proposition}
\newtheorem{corollary}[theorem]{Corollary}
\newtheorem{facts}[theorem]{Facts}
\newtheorem{props}[theorem]{Properties}
\newtheorem*{thmA}{Theorem A}
\newtheorem{ex}[theorem]{Example}
\theoremstyle{definition}
\newtheorem{definition}[theorem]{Definition}
\newtheorem{remark}[theorem]{Remark}
\newtheorem*{defna}{Definition}

\newcommand{\notes} {\noindent \textbf{Notes.  }}
\newcommand{\note} {\noindent \textbf{Note.  }}
\newcommand{\defn} {\noindent \textbf{Definition.  }}
\newcommand{\defns} {\noindent \textbf{Definitions.  }}
\newcommand{\x}{{\bf x}}
\newcommand{\z}{{\bf z}}
\newcommand{\B}{{\bf b}}
\newcommand{\V}{{\bf v}}
\newcommand{\T}{\mathbb{T}}
\newcommand{\Z}{\mathbb{Z}}
\newcommand{\Hp}{\mathbb{H}}
\newcommand{\D}{\mathbb{D}}
\newcommand{\R}{\mathbb{R}}
\newcommand{\N}{\mathbb{N}}
\renewcommand{\B}{\mathbb{B}}
\newcommand{\C}{\mathbb{C}}
\newcommand{\ft}{\widetilde{f}}
\newcommand{\dt}{{\mathrm{det }\;}}
 \newcommand{\adj}{{\mathrm{adj}\;}}
 \newcommand{\0}{{\bf O}}
 \newcommand{\av}{\arrowvert}
 \newcommand{\zbar}{\overline{z}}
 \newcommand{\xbar}{\overline{X}}
 \newcommand{\htt}{\widetilde{h}}
\newcommand{\ty}{\mathcal{T}}
\renewcommand\Re{\operatorname{Re}}
\renewcommand\Im{\operatorname{Im}}
\newcommand{\tr}{\operatorname{Tr}}

\newcommand{\ds}{\displaystyle}
\numberwithin{equation}{section}

\renewcommand{\theenumi}{(\roman{enumi})}
\renewcommand{\labelenumi}{\theenumi}

\title{On Infinitesimal Strebel Points}
\author{Alastair Fletcher}
\address{Department of Mathematical Sciences, Northern Illinois University, DeKalb, IL 60115-2888. USA}
\email{fletcher@math.niu.edu}

\thanks{This work was supported by a grant from the Simons Foundation (\#352034, Alastair Fletcher).}
\maketitle

\begin{abstract}
In this paper, we prove that if $X$ is a Riemann surface of infinite analytic type and $[\mu ]_T$ is any element of Teichm\"uller space, then there exists $\mu_1 \in [\mu]_T$ so that $\mu_1$ is an infinitesimal Strebel point.

\end{abstract}

\section{Introduction}

If $X$ is a hyperbolic Riemann surface and $M(X)$ denotes the set of Beltrami differentials on $X$, then Teichm\"uller space $T(X)$ is defined by the set of equivalence classes in $M(X)$ under the Teichm\"uller equivalence relation.
Given $\mu \in M(X)$, we denote by $[\mu]_T \in T(X)$ the corresponding point of Teichm\"uller space.
In the case where $X$ is a closed surface, every $[\mu ]_T$ different from the basepoint $[0]_T$ contains a uniquely extremal representative of the form 
\begin{equation}
\label{eq:mu}
\mu = k|\varphi|/\varphi
\end{equation} 
for $0< k <1$ and $\varphi \in A^1(X)$, that is, $\varphi$ is an integrable holomorphic quadratic differential on $X$.

In the setting of infinite type surfaces, not every Teichm\"uller class has a uniquely extremal representative. The first such example to be constructed was the Strebel chimney \cite{S}. If $X = \{ z\in \C : \Im (z) <0 \text{ or } |\Re(z)|<1 \}$ and $K>1$, then for every $L\in [1/K,K]$
\[ f_L(x+iy) = \left \{ \begin{array}{cc} x+iKy, & y\geq 0, |x|<1 \\ x+iLy, & y < 0 \end{array} \right .\]
is an extremal representative in its Teichm\"uller class $[\mu_{f_1}]_T$. Here $\mu_f \in M(X)$ denotes the complex dilatation of a quasiconformal mapping $f$.

Strebel points, introduced in \cite{EL,L}, are those points in an infinite dimensional Teichm\"uller space where the behaviour is, in some sense, tame. In particular, Strebel's frame mapping criterion \cite{G} states that every Strebel point can be represented by a uniquely extremal Beltrami differential of the form \eqref{eq:mu}. It is known that the set of Strebel points is both open and dense in Teichm\"uller space.

We may also consider the infinitesimal class of $\mu \in M(X)$, namely those $\nu \in M(X)$ which induce the same linear functional $\Lambda _{\nu}(\varphi) = \int_X \nu \varphi$ on $A^1(X)$ as $\mu$. The set of equivalence classes is denoted by $B(X)$ and a particular class is denoted $[\mu]_B$. We can consider extremality of representatives of an equivalence class in $B(X)$. If $\mu$ is extremal in its class in $B(X)$, it is called infinitesimally extremal.
There is an analogous notion of an infinitesimal Strebel point, and the infinitesimal version of Strebel's frame mapping criterion (see \cite{EG}) states that for such a point, it can again be represented by a uniquely infinitesimally extremal Beltrami differential of the form \eqref{eq:mu}.

Typically $[\mu]_T $ and $[\mu]_B$ do not coincide but they share similar properties. In particular, $\mu$ is extremal if and only if it is infinitesimally extremal, see \cite{H,K,RS}. Further, $\mu$ is uniquely extremal if and only if it is uniquely infinitesimally extremal, see \cite{BLMM}. A point to note here is that not every uniquely extremal Beltrami differential has the form \eqref{eq:mu}.

In \cite{HS}, it was shown that if $||\mu ||_{\infty}$ is small then there exists $\nu_1 \in [\mu ]_T$ such that $[\nu_1]_B$ is an infinitesimal Strebel point and $\nu_2 \in [\mu ]_B$ such that $[\nu_2]_B$ is a Strebel point. Further, in \cite{HL} it was proved that for any $\mu \in M(\D)$ there exists $\nu \in [\mu]_T$ such that $\nu$ is an infinitesimal Strebel point. In this paper, we extend this latter result to all surfaces of infinite type.

\begin{theorem}
\label{thm:main}
Let $X$ be a Riemann surface of infinite analytic type. Given any Beltrami differential $\mu \in M(X)$, there exists $\mu_1 \in [\mu]_{T}$ such that $\mu_1$ is an infinitesimal Strebel point.
\end{theorem}

Before proving this theorem, we first recall some material from Teichm\"uller theory, see for example \cite{FM,GL} for more details.

{\it Acknowledgements:} The author would like to thank the referee for pointing out a gap in an earlier version of this paper.

\section{Preliminaries}

\subsection{Beltrami differentials}
Let $X$ be a Riemann surface of infinite analytic type, that is, we require the genus of $X$ or the number of punctures to be infinite. Denote by $M(X)$ the set of Beltrami differentials on $X$, that is, $(-1,1)$-differential forms given in local coordinates by 
\[ \mu(z) \frac{ d\overline{z}}{dz},\]
and where $||\mu||_{\infty} \leq k$ for some $k<1$.

There are two equivalence relations we may impose on $M(X)$. The Teichm\"uller equivalence relation is defined as follows. Since $X$ is a hyperbolic surface, it can be realized as $\D / G$ for some group $G$ of covering transformations.
We may lift $\mu,\nu \in M(X)$ to elements $\widetilde{\mu}, \widetilde{\nu}$ in the unit ball of $L^{\infty} (\D) $ which are $G$-invariant in the sense that
\[ (\widetilde{\mu} \circ g) \frac{ \overline{g'}}{g'} = \widetilde{\mu},\]
for all $g\in G$. Via the Measurable Riemann Mapping Theorem, there is a unique quasiconformal mapping $f^{\widetilde{\mu}} : \D \to \D$ which solves the Beltrami differential equation $f_{\zbar} = \widetilde{\mu} f_z$, extends continuously to the boundary $\partial \D$ and fixes $1,-1,i$. We then say that $\mu \sim_T \nu$ if and only if
\begin{equation}
\label{eq:class} 
f^{\widetilde{\mu}} |_{\partial \D} = f^{\widetilde{\nu}} |_{\partial \D},
\end{equation}
and denote an equivalence class by $[\mu]_T$. The set of equivalence classes is the Teichm\"uller space $T(X)$ of $X$.

The second equivalence class is defined as follows. Let $A^1(X)$ be the Bergman space of integrable holomorphic quadratic differentials on $X$, that is, $(2,0)$-differential forms give in local coordinates by
\[ \varphi (z) \: dz^2,\]
where $\varphi$ is holomorphic and satisfies
\[ ||\varphi||_1 := \int_X | \varphi| < \infty .\]
We remark that every such $\varphi$ lifts to a holomorphic function $\widetilde{\varphi}$ on $\D$ satisfying $(\widetilde{\varphi} \circ g)(g')^2 = \widetilde{\varphi}$ for all $g\in G$.
It is well-known that the cotangent space at the basepoint of Teichm\"uller space, $[0]_T$, is isomorphic to $A^1(X)$. Every $\mu \in M(X)$ induces a linear functional on $A^1(X)$ defined by
\[ \Lambda _{\mu} (\varphi) = \int _X \mu \varphi.\]
Note that the expression $\mu \varphi$ in local coordinates is $\mu(z) \varphi(z) |dz|^2$ and so can be integrated. If $\Lambda _{\mu} = \Lambda_{\nu}$ on $A^1(X)$, then we say that $\mu$ and $\nu $ are infinitesimally equivalent and write $\mu \sim_B \nu$. An equivalence class is denoted by $[\mu]_B$ and the set of equivalence classes is denoted by $B(X)$, and is isomorphic to $A^1(X)^*$. The norm of $\Lambda_{\mu}$ is
\[ ||\Lambda_{\mu} || = \sup \left \{ \left | \int_X \mu \varphi \right | : ||\varphi||_1=1 \right \}.\]
If $\mu$ is extremal, then $||\Lambda _{\mu}|| = ||\mu ||_{\infty}$.

In either equivalence class, we have the notion of extremality and unique extremality. If $\sim$ is either $\sim_T$ or $\sim_B$, then we say that $\mu$ is extremal if $||\mu||_{\infty} \leq ||\nu ||_{\infty}$ for all $\nu \sim \mu$, and we say that $\mu$ is uniquely extremal if the inequality is strict for all $\nu \sim \mu$ and $\nu \neq \mu$. We will use the language {\it extremal} for $\sim_T$ and {\it infinitesimally extremal} for $\sim_B$. From compactness considerations, there is always an extremal and infinitesimal extremal representative in each class.

\subsection{Extremal distortion}
The notation we use in this section is slightly non-standard but is intended to unify the ideas between the two equivalence classes.
Writing $k(\mu) = || \mu ||_{\infty}$ for $\mu \in M(X)$, the corresponding extremal version is
\[ k_0( [\mu]_T) = \inf \{ k(\nu) : \nu \sim_T \mu \}  \]
and the infinitesimally extremal version is
\[ k_1( [\mu]_B) = \inf \{ k(\nu) : \nu \sim _B \mu \}.\]
It is known that
\begin{equation}
\label{eq:B} 
k_1([\mu]_B) = ||\Lambda_{\mu} || =  \sup_{ ||\varphi||_1=1} \Re \int_X \mu \varphi.
\end{equation}
The boundary dilatation of $\mu \in M(X)$ is 
\[ h( \mu) = \inf \{ || \mu |_{X\setminus E} ||_{\infty} : E \text{ is a compact subset of } X \}.\]
The extremal version is 
\[ h_0( [\mu ]_T) = \inf \{h(\nu ) : \nu \sim_T \mu \} \]
and the infinitesimally extremal version is 
\[ h_1( [\mu ]_B) = \inf \{ h(\nu ) : \nu \sim_B \mu \}.\]
It is clear that $h_0([\mu]_T) \leq k_0([\mu]_T)$ and $h_1([\mu]_B) \leq k_1([\mu]_B)$. In the first case, if $h_0([\mu]_T) < k_0([\mu]_T)$ then $[\mu]_T \in T(X)$ is called a {\it Strebel point}, and otherwise a non-Strebel point. In the second case, if $h_1([\mu]_B) < k_1([\mu]_B)$, then $[\mu]_B \in B(X)$ is called an {\it infinitesimal Strebel point}, and otherwise an infinitesimal non-Strebel point.

The case $\mu = 0$ gives a non-Strebel point $[0]_T \in T(X)$ and an infinitesimal non-Strebel point $[0]_B \in B(X)$.

\subsection{Quasiconformal gluing}

We will require the following result on quasiconformal gluing. 

\begin{theorem}[\cite{JQ}, Theorem 2]
\label{thm:gluing}
Let $f_1$ and $f_2$ be two $K$-quasiconformal mappings defined on disjoint simply connected subdomains $\Omega_1$ and $\Omega_2$ of $\overline{\C}$ respectively and $f_1(\Omega_1) \cap f_2( \Omega_2)  =\emptyset$. Then for any two Jordan domains $D_1$ and $D_2$ with $\overline{D_1} \subset \Omega_1$ and $\overline{D_2} \subset \Omega_2$, there exists a quasiconformal mapping $g$ of $\overline{\C}$ such that
\[ g|_{D_1} = f_1|_{D_1} \: \text{ and } \: g|_{D_2}  = f_2|_{D_2}.\]
\end{theorem}

\section{Proof of Theorem \ref{thm:main}}

In \cite{HL}, where the case $X = \D$ is considered, a particular choice of $\varphi = 1/\pi \in A^1(\D)$ is chosen to simplify estimates. In our situation, there isn't an obvious choice of $\varphi$. To deal with this, we will choose any $\varphi$ of norm $1$, decompose our surface $X$ into three subsets and modify a given $\mu$ on some of these subsets. The modified Beltrami differential will remain in the same Teichm\"uller class but will be an infinitesimal Strebel point.

Fix a Riemann surface $X$ of infinite analytic type and $\varphi \in A^1(X)$ with $||\varphi ||_1=1$. Let $\mu \in M(X)$ with $||\mu||_{\infty}= k <1$ and $\epsilon >0$.

{\bf Step 1:} find an appropriate decomposition of $X$. Realize $X$ as a union of pairs of pants, that is, topological three-holed spheres where some of the holes are allowed to be points, see for example \cite{AR}. Label the pants $P_1,P_2,\ldots$ where for convenience we may assume $\bigcup_{i=1}^n P_i$ is connected for each $n$. Choose $N$ large enough that if $G_1 = \bigcup_{i=N+1}^{\infty} P_i$, then 
\[ \int _{G_1} |\varphi | <\frac{\epsilon}{2}.\]

Each $P_i$ is obtained by gluing two hypergolic geodesic hexagons together along three pairs of sides called seams. The three boundary components of $P_i$ are called cuffs. For $i=1,\ldots,N$, let $U_i$ be an open neighbourhood of the cuffs and seams in $P_i$ so that $P_i \setminus U_i$ consists of two components, one arising from each hexagon, and so that
\[ \int_{U_i} |\varphi| <\frac{\epsilon}{4N}.\]
Note there is no issue if a boundary component of $P_i$ reduces to a point. In this case, $U_i$ wil contain a neighbourhood of this puncture.

We next take closed neighbourhoods of each $U_i$, say $F_i$, which will consist of two components, each a topological annulus so that $E_i : = P_i \setminus (U_i \cup F_i)$ consists of two topological disks. We further require that
\[ \int _{F_i} |\varphi| <\frac{ \epsilon}{4N}.\]

Our disjoint decomposition of $X$ is then defined by $G = G_1 \cup \bigcup_{i=1}^N U_i$, $F = \bigcup_{i=1}^N F_i$ and $E = \bigcup _{i=1}^N E_i$. By construction, $E$ and $G$ are open, $F$ is closed and
\begin{equation}
\label{eq:estimates} 
\int_G |\varphi| = \int_{G_1} |\varphi| + \sum_{i=1}^N \int_{U_i} |\varphi| <\frac{3\epsilon}{4} ,\: 
\int _F |\varphi| = \sum_{i=1}^N \int_{F_i} |\varphi| < \frac{\epsilon}{4}, \:
\int_E |\varphi| >1-\epsilon.
\end{equation}

{\bf Step 2:} Given a Beltrami differential $\mu$ in $M(X)$, we modify it on $E$ and $F$. Let $f:X \to Y$ be a quasiconformal map with complex dilatation $\mu$. Denote by $\pi_1,\pi_2$ the projections from $\D$ to $X,Y$ respectively.
We may lift $f$ to a quasiconformal map $\widetilde{f} :\D \to \D$ which satisfies $\pi_2 \circ \widetilde{f} = f\circ \pi_1$. Since $\widetilde{f}$ extends continuously to $\partial \D$, we may extend $\widetilde{f}$ to a quasiconformal map $\widetilde{f} :\overline{\C} \to \overline{\C}$ via reflection in $\partial \D$.

Given a pair of pants $P_i$ with $i\in \{1,\ldots,N\}$, consider one of the two hexagons $\sigma$ that are glued to give $P_i$. 
Consider a component $H$ of $\pi_1^{-1}(\sigma)$, which is a geodesic hexagon in $\D$.
Let $H' = \widetilde{f}(H) \subset \D$.
If we intersect $H$ with $\pi_1^{-1}(E)$ and $\pi_1^{-1}(F)$ then we obtain an open topological disk $D$ and a closed topological annulus $A$ respectively. Let $\Omega_1 \subset \overline{\C}$ be an open neighbourhood of $\overline{\C} \setminus (A \cup D)$ and let $\Omega_2 \subset \overline{\C}$ be an open neighbourhood of $\overline{D}$ so that $\overline{\Omega_1} \cap \overline{\Omega_2} = \emptyset$ and both $\Omega_1$ and $\Omega_2$ are simply connected in $\overline{\C}$. Note that $\partial \Omega_1$ and $\partial \Omega_2$ are both contained in the interior of $A$.

By the Measurable Riemann Mapping Theorem, we may find a quasiconformal map $g:\Omega_2 \to f(\Omega_2)$ with complex dilatation
\[ \mu_1 = \frac{k_1\overline{\widetilde{\varphi}}}{|\widetilde{\varphi} |},\]
where $k_1$ is to be determined and $\widetilde{\varphi}$ is a lift to $\D$ of the Beltrami differential $\varphi$ on $X$.
We may then apply Theorem \ref{thm:gluing} with $\Omega_1,\Omega_2$ as above, $f_1 = \widetilde{f}$, $f_2 = g$ and $D_1 = \overline{\C} \setminus (A \cup D)$, $D_2 = D$.
Hence there is a quasiconformal map $h:\overline{\C} \to \overline{\C}$ satisfying
\[ h|_D = g, \: h|_{\overline{\C} \setminus (A \cup D)} = \widetilde{f}.\]
We then replace $\widetilde{f}$ on $H$ by this new map $h$.

Repeating this construction on the lift of each hexagon that makes up $P_1,\ldots,P_N$, we end up with a quasiconformal map defined on a fundamental region of $X$ in $\D$ with image a fundamental region for $Y$ in $\D$. We can then obtain a quasiconformal map $\widetilde{f_1} :\D \to \D$ by propagating to other fundamental regions in the domain and range via the covering groups.
Projecting back to $X$, we obtain a quasiconformal map $f_1$ on $X$ which agrees with $f$ outside a compact set and on a neighbourhood of the seams and cuffs, and agrees with a quasiconformal map with complex dilatation $k_1\overline{\varphi} / | \varphi |$ on a large subset. We will denote by $\mu_1$ the complex dilatation of $f_1$. 

{\bf Step 3:} $\mu_1$ is in the same Teichm\"uller class as $\mu$. Recall that $\widetilde{f_1}$ is a lift of $f_1$ to $\D$. Since $\widetilde{f_1}$ extends continuously to the boundary, we just need to show that the boundary map of $\widetilde{f_1}$ agrees with that of $\widetilde{f}$. To see this, let $w_0 \in \partial \D$. Since the components of $\pi_1^{-1}(X\setminus G)$ in $\D$ are pairwise disjoint and have uniformly bounded hyperbolic diameter,
we can find a sequence of points $z_n$ contained in $\pi_1^{-1}(G)$ which converges to $w_0$ in the Euclidean metric. Since $\widetilde{f_1}$ agrees with $\widetilde{f}$ on $\pi_1^{-1}(G)$ and both maps extend continuously to $\partial \D$, we see that $\widetilde{f_1}(w_0) = \widetilde{f}(w_0)$. Since $w_0$ is arbitrary, by \eqref{eq:class} it follows that these two maps determine the same two points of Teichm\"uller space.

{\bf Step 4:} $\mu_1$ is infinitesimally extremal for an appropriate choice of $k_1$. Choose $\epsilon>0$ small enough and let 
\begin{equation}
\label{eq:2}
k_1 > \frac{(1+\epsilon) k + \epsilon }{1-\epsilon}.
\end{equation}
Note that since the function $p(x) = ((1+x)k+x)/(1-x)$ satisfies $p(0) = k$ and $p'(0) >0$, this can be achieved.

We have by \eqref{eq:B}
\begin{align*}
|| \Lambda _{\mu_1} || & = \sup \left \{ \left | \Re  \int _X \mu_1 \psi \right | : \psi \in A^1(X), ||\psi||_1=1 \right \} \\
& \geq \left | \Re \int _X \mu_1 \varphi \right | \\
& \geq  \left | \Re \int _E \mu_1 \varphi  +\Re \int _F \mu_1 \varphi +\Re \int _G \mu_1 \varphi \right   |\\
& \geq \left | \Re \int _E \mu_1 \varphi \right |  - \left | \Re \int _F \mu_1 \varphi \right | - \left | \Re \int _G \mu_1 \varphi \right  |.
\end{align*}
By \eqref{eq:estimates}, we have
\[ \Re \int_E \mu_1 \varphi  = \int_E \mu_1 \varphi  = \int _E k_1 |\varphi | > k_1(1-\epsilon),\]
and
\[ \left | \Re \int _F \mu_1 \varphi \right | \leq ||\mu_1 ||_{\infty} \int_F | \varphi | < \epsilon.\]
The contribution on $G$ is, again by \eqref{eq:estimates}
\[ \left | \Re \int _G \mu_1 \varphi \right | = \left | \Re \int_G \mu \varphi  \right | \leq ||\mu ||_{\infty} \int _G |\varphi | < k \epsilon.\]
Combining these estimates, we obtain
\[ || \Lambda _{\mu_1} || > k_1(1-\epsilon )  - (k+1)\epsilon.  \]
Finally, by \eqref{eq:2}, we obtain 
\[ k_1( [\mu_1]_B) =  || \Lambda _{\mu_1} || > k_1(1-\epsilon )  - (k+1)\epsilon >  k > || \mu ||_{\infty} \geq || \mu_1 |_G ||_{\infty} \geq h_1([\mu_1]_B),\]
and so $\mu_1$ is an infinitesimal Strebel point.

\end{document}